\documentclass[11pt]{article}

\usepackage{amscd, amsmath, amssymb, amsthm, fancyhdr, graphicx, epsfig, url, mathrsfs}
\usepackage{graphicx}
\usepackage{epigraph}
\usepackage{tikz-cd}
\usepackage[backref=page]{hyperref}
\renewcommand*{\backrefalt}[4]{%
	\ifcase #1 (Not cited.)%
	\or        (Cited on page~#2.)%
	\else      (Cited on pages~#2.)%
	\fi}

\hypersetup{
	colorlinks   = true,
	citecolor    = magenta,
	linkcolor    = blue,
	urlcolor     = magenta	
}

\numberwithin{equation}{section}

\newcommand{\version}{version 1.1.1,\ \ October 28, 2025}

\makeatletter
\def\x@arrow{\DOTSB\Relbar}
\def\xlongrightarrowfill@{\arrowfill@\relbar\relbar\longrightarrow}
\newcommand{\xlongrightarrow}[2][]{%
        \ext@arrow 0099\xlongrightarrowfill@{#1}{#2}}
\makeatother

\newcommand{\dra}{\dashrightarrow}

\def\eqref#1{(\ref{#1})}

\newcommand{\arrow}{{\:\longrightarrow\:}}

\def\1{\sqrt{-1}\:}

\newcommand{\cntrct}                
{\hspace{2pt}\raisebox{1pt}{\text{$\lrcorner$}}\hspace{2pt}}

\renewcommand{\tilde}{\widetilde}
\renewcommand{\bar}{\overline}
\renewcommand{\phi}{\varphi}
\renewcommand{\epsilon}{\varepsilon}

\renewcommand{\max}{{\rm max}}

\newcommand{\Teich}{\operatorname{\sf Teich}}

\newcommand{\MBM}{\operatorname{MBM}}

\newcommand{\Hom}{\operatorname{Hom}}

\newcommand{\Pic}{\operatorname{Pic}}
\newcommand{\Pos}{\mathcal{C}}

\newcommand{\codim}{\operatorname{codim}}

\newcommand{\rk}{\operatorname{rk}}

\newcommand{\emrp}{\operatorname{End}}

\newcommand{\BK}{\mathcal{BK}}

\newcommand{\bbZ}{\mathbb{Z}}
\newcommand{\bbQ}{\mathbb{Q}}
\newcommand{\bbR}{\mathbb{R}}
\newcommand{\bbC}{\mathbb{C}}

\newcommand{\bbP}{\mathbb{P}}

\newcommand{\CC}{\mathcal{C}}
\newcommand{\DD}{\mathcal{D}}

\newcommand{\KK}{\mathcal{K}}
\newcommand{\LL}{\mathcal{L}}
\newcommand{\MM}{\mathcal{M}}
\newcommand{\NN}{\mathcal{N}}
\newcommand{\OO}{\mathcal{O}}

\renewcommand{\ge}{\geqslant}
\renewcommand{\le}{\leqslant}

\newcommand{\st}{\enskip |\enskip}

\newcommand{\cnv}{\,\lrcorner\,}

\newcommand{\wdg}{\wedge}

\newcommand{\hrarr}{\hookrightarrow}

\newcounter{Mycounter}[section]
\newcounter{lemma}[section]
\renewcommand{\thelemma}{{Lemma \thesection.\arabic{lemma}}}
\newcommand{\lemma}{%
    \setcounter{lemma}{\value{Mycounter}}
    \refstepcounter{lemma}
    \stepcounter{Mycounter}
    {\noindent \bf \thelemma:\ }}

\newcounter{claim}[section]
\newcounter{sublemma}[section]
\newcounter{corollary}[section]
\renewcommand{\thecorollary}{{Corollary \thesection.\arabic{corollary}}}
\newcommand{\corollary}{%
    \setcounter{corollary}{\value{Mycounter}}
    \refstepcounter{corollary}
    \stepcounter{Mycounter}
    {\noindent \bf \thecorollary:\ }}

\newcounter{theorem}[section]
\renewcommand{\thetheorem}{{Theorem \thesection.\arabic{theorem}}}
\newcommand{\theorem}{%
    \setcounter{theorem}{\value{Mycounter}}
    \refstepcounter{theorem}
    \stepcounter{Mycounter}
    {\noindent \bf \thetheorem:\ }}

\newcounter{conjecture}[section]
\renewcommand{\theconjecture}{{Conjecture \thesection.\arabic{conjecture}}}
\newcommand{\conjecture}{%
    \setcounter{conjecture}{\value{Mycounter}}
    \refstepcounter{conjecture}
    \stepcounter{Mycounter}
    {\noindent \bf \theconjecture:\ }}

\newcounter{proposition}[section]
\renewcommand{\theproposition}
      {{Proposition \thesection.\arabic{proposition}}}
\newcommand{\proposition}{%
    \setcounter{proposition}{\value{Mycounter}}
    \refstepcounter{proposition}
    \stepcounter{Mycounter}
    {\noindent \bf \theproposition:\ }}

\newcounter{definition}[section]
\renewcommand{\thedefinition}
      {{Definition~\thesection.\arabic{definition}}}
\newcommand{\definition}{%
    \setcounter{definition}{\value{Mycounter}}
    \refstepcounter{definition}
    \stepcounter{Mycounter}
    {\noindent \bf \thedefinition:\ }}

\newcounter{example}[section]
\newcounter{remark}[section]
\newcounter{problem}[section]
\newcounter{question}[section]
\makeatletter

\@addtoreset{equation}{section} \@addtoreset{footnote}{section}
\makeatother

\newtheorem{defn}{Definition}[section]

{\theoremstyle{remark}
  \newtheorem{rem}[defn]{Remark}
  
}

\begin{document}

\begin{center}
{\LARGE\bf
Birational geometry of hyperk\"ahler\\[3mm] manifolds and the Hu-Yau conjecture
}

\bigskip

{Ekaterina Amerik,
Andrey Soldatenkov\footnote{Partially supported by 
FAPESP grant 2024/23819-0}, 
Misha Verbitsky\footnote{Partially supported by 
FAPERJ grant SEI-260003/000410/2023 and CNPq - Process 310952/2021-2. \\


 }}

\end{center}

{\small\hspace{.15\linewidth}\begin{minipage}[t]{0.8\linewidth}
{\bf Abstract:} 
Wierzba and
Wi\'sniewski proved that for $\dim_\bbC M=4$, every
bimeromorphic map $M \dra M'$ of hyperk\"ahler manifolds
is represented as a composition of Mukai flops.
Hu and Yau conjectured that this result can
be generalized to arbitrary dimension. They
defined ``Mukai's elementary transformation'' as the
blow-up of a subvariety ruled by complex projective spaces, composed with
the contraction of the ruling. Hu and Yau conjectured that 
any bimeromorphic map $M \dra M'$ can be decomposed into a sequence of Mukai's
elementary transformations, after possibly removing
subvarieties of codimension greater than $2$.
We prove this conjecture for compact hyperk\"ahler manifolds
of maximal holonomy by decomposing any bimeromorphic map into a
composition of wall-crossing flops associated with MBM contractions.
\end{minipage}
}

\tableofcontents


\section{Introduction}


Let $M$ and $M'$ be two bimeromorphic compact hyperk\"ahler manifolds.
In \cite{_HuYau:AdvTeor_} Hu and Yau have studied the
most basic type of bimeromorphic maps between $M$ and $M'$
known as Mukai flops or Mukai's elementary transformations.
Hu and Yau conjectured that every bimeromorphic map between
$M$ and $M'$ may to a certain extent be described as a sequence
of Mukai flops. More precisely, given a bimeromorphic map
$\varphi\colon M \dra M'$, it was conjectured
that after removing some subvarieties of codimension
greater than two from $M$ and $M'$ the map $\varphi$ decomposes into
a finite composition of Mukai's elementary transformations.

The work by Wisniewski and Wierzba \cite{WW}, along with some observations in \cite{BHL}, actually proves a stronger result in dimension four. In the case of fourfolds one does not need
to remove any subvarieties from $M$ and $M'$ to decompose $\varphi$
into a sequence of flops. The central result of \cite[Theorem 1.1]{WW} is the
description of the contraction loci in holomorphic symplectic fourfolds,
claiming that the contraction locus is always a disjoint union of
complex projective planes.

In the present paper we address the Hu--Yau conjecture in arbitrary
dimension. Our aim is to give a precise statement of the conjecture,
drawing attention to the possible ambiguity in the formulation.
The ambiguity originates from the need to remove subvarieties
of codimension three and higher from the original hyperk\"ahler
manifolds. To make things precise, we introduce in Section \ref{sec_HY} the notion
of a Hu--Yau transformation: this is a Mukai's elementary transformation up to codimension three or higher. Our main result, \ref{thm_main}, claims that any bimeromorphic map $\varphi$
as above is the composition of a finite number of Hu--Yau transformations.

Let us make a few comments on the comparison of our proof with that of \cite{WW}.
Our proof relies on \cite[Theorem 1.1]{WW} in an essential way. On the other hand,
our methods give a new proof of \cite[Theorem 1.2]{WW} in the case when the
varieties are compact hyperk\"ahler. The novelty here is a different method
of decomposing bimeromorphic maps into compositions of flops. Whereas loc. cit.
relies on the general arguments of the minimal model program providing the termination
of log-flips, our proof uses completely different methods based on \cite[Theorem 3.17]{_AV:Orbits_}.
The latter theorem (see also \cite[Proposition 17]{HT-Moving} provides a locally finite wall-and-chamber structure for the
positive cone of an arbitrary compact hyperk\"ahler manifold (under mild conditions on its second Betti number),
and therefore the possibility to decompose an arbitrary bimeromorphic map into a finite sequence of wall-crossing transformations,
see Section \ref{sec_WC} for details. Unlike \cite[Theorem 1.2]{WW}, our proof works in the non-projective
setting, therefore generalizing the main results of loc. cit. to the case of non-projective
hyperk\"ahler fourfolds. 

In Section \ref{sec_compositions} we provide an example to justify our reformulation of the
Hu--Yau conjecture. The example deals with the geometry of intersections of the centers of Mukai's elementary transformations
in codimension two, and to the dynamical behaviour of the composition of such transformations.


\section{Main results}


In this section, we give the key definitions, in order to state the main
result of this paper, \ref{thm_main}. 


\subsection{Mukai flops} \label{_Mukai_flop_Subsection_}


Given a complex $(k + 1)$-dimensional vector space $V$,
let $M_{2k} = \mathrm{Tot}(\Omega^1_{\bbP(V)})$ and $M'_{2k} = \mathrm{Tot}(\Omega^1_{\bbP(V^*)})$
be the total spaces of the cotangent bundles of $\bbP (V)$ and of its dual.
It is well-known that the zero sections of the bundles $M_{2k}$ and $M'_{2k}$ can be contracted,
and there is an isomorphism between the resulting varieties that gives rise to a birational map
between $M_{2k}$ and $M'_{2k}$ called the Mukai flop.

Let us briefly recall some details of this construction. Starting from the Euler sequence on $\bbP(V)$:
$$
0\to \OO_{\bbP(V)} \to  V \otimes \OO_{\bbP(V)}(1) \to T_{\bbP(V)} \to 0,
$$
we pass to its dual and use the canonical map $\OO_{\bbP(V)}(-1) \hrarr V\otimes \OO_{\bbP(V)}$ to get the embedding of vector bundles
$$
\Omega^1_{\bbP(V)} \hrarr \emrp(V)\otimes \OO_{\bbP(V)}.
$$
Given a point $[u]\in \bbP(V)$ recall that $\Omega^1_{\bbP(V),[u]} \simeq \Hom(V/\bbC u, \bbC u)$.
We see that $M_{2k}$ can be described as the following closed subset of $\bbP(V)\times \emrp(V)$:
$$
M_{2k} = \{([u], A)\in \bbP(V)\times\emrp(V)\st \mathrm{im}(A)\subset \bbC u,\, A^2 = 0 \}.
$$
The projection $\pi_{2k}\colon M_{2k}\to \emrp(V)$ contracts the zero section $$Z_{2k} = \{([u], A)\in M_{2k}\st A = 0\}$$ to the origin.
The image $N_{2k} = \pi_{2k}(M_{2k})$ is the closure of the minimal nilpotent $\mathrm{GL}(V)$-orbit in $\emrp(V)$
consisting of square-zero endomorphisms of rank at most one. The symplectic form on $M_{2k}$ is the pullback of the Kirillov--Kostant form on $N_{2k}$.
Analogously,
$$
M_{2k}' = \{([\alpha], B)\in \bbP(V^*)\times\emrp(V^*)\st \mathrm{im}(B)\subset \bbC \alpha,\, B^2 = 0 \}.
$$
There is the analogous projection $\pi_{2k}'\colon M_{2k}'\to \emrp(V^*)$ onto the closed affine cone $N_{2k}'\subset \emrp(V^*)$.
Identifying $\emrp(V)$ with $\emrp(V^*)$ by taking adjoints of the endomorphisms gives an isomorphism between $N_{2k}$ and $N_{2k}'$.

\hfill

\definition\label{defn_Mukai_flop}
The morphism $\mu\colon M_{2k}\setminus Z_{2k}\to M_{2k}'$ defined by $\mu([u], A) = ([\ker(A)], A^*)$
gives a birational map
\begin{equation}\label{eqn_Mukai_flop}
\mu\colon M_{2k} \dra M_{2k}'
\end{equation}
called the {\bf Mukai flop} in codimension $k$. It is clear from the
above description that $\mu = (\pi_{2k}')^{-1}\circ\pi_{2k}$, where $\pi_{2k}\colon M_{2k}\to N_{2k}$ and $\pi_{2k}'\colon M_{2k}'\to N_{2k}'\simeq N_{2k}$
are the contractions described above.

\hfill

Given a hyperk\"ahler fourfold $M$ and a projective plane $S\simeq \bbC P^2\subset M$,
an open analytic neighborhood of $S$ in $M$ is isomorphic to a neighbourhood of the zero section $Z_4$ in $M_4$,
see \cite{WW} and references therein, hence the Mukai flop construction can be performed
on any hyperk\"ahler fourfold. In \cite{WW}, J. Wierzba and J. Wi\'sniewski prove that a small contraction
of a quasiprojective holomorphic symplectic fourfold contracts a disjoint union of $\bbP^2$'s, and it follows that
any bimeromorphic map between hyperk\"ahler manifolds of complex dimension four can be
decomposed into a composition of Mukai flops (\cite{BHL}).

More generally, for hyperk\"ahler manifolds of arbitrary dimension,
we would like to describe the building blocks of arbitrary bimeromorphic maps.
The local model of our bimeromorphic transformations will be given by taking the product of
a Mukai flop with a polydisc (possibly up to codimension three or higher, see below).

More precisely, let
$U_{2k}\subset N_{2k}$ be an open analytic neighborhood of the vertex $0\in N_{2k}$,
where $N_{2k}$ is the affine cone from \ref{defn_Mukai_flop}. Using the morphisms
$\pi_{2k}$ and $\pi'_{2k}$ from \ref{defn_Mukai_flop}, we define the local model
of a bimeromorphic transformation of a $d$-dimensional complex manifold by the following
diagram:
\begin{equation}\label{eqn_local_model}
\begin{tikzcd}[column sep = small]
  \pi^{-1}_{2k}(U_{2k})\times\Delta^{d-2k} \arrow{dr}\arrow[dashed]{rr}{\mu\times \mathrm{id}} & & (\pi'_{2k})^{-1}(U_{2k})\times\Delta^{d-2k} \arrow{dl} \\
  & U_{2k}\times\Delta^{d-2k} &
\end{tikzcd}
\end{equation}
In the above diagram, the diagonal arrows define contractions with the exceptional
sets isomorphic to $\bbP^k$-bundles over the polydisc $\Delta^{d-2k}$. We introduce
the following notion, using the terminology of \cite{_HuYau:AdvTeor_}.

\hfill

\definition\label{defn_Mukai_flop2}
Let $\phi\colon X \dra X'$ be a bimeromorphic map between arbitrary
complex manifolds of complex dimension $d$. We will say that $\phi$ is {\bf Mukai's elementary transformation}
in codimension $k$, if there exists a normal complex-analytic variety $Y$, two closed non-empty complex
submanifolds $E\subset X$, $E'\subset X'$ and two proper bimeromorphic morphisms
$\pi\colon X\to Y$, $\pi'\colon X'\to Y$ that satisfy the following conditions.
\begin{enumerate}
\item $\phi = (\pi')^{-1}\circ \pi$;
\item $\pi$ and $\pi'$ are contractions of $E$ and $E'$ onto the singular locus $Y^s\subset Y$;
\item for any point $y\in Y^s$ there exists an open neighborhood $U\subset Y$
of $y$ isomorphic to $U_{2k}\times \Delta^{d-2k}$, and over $U$ the map $\phi$
is of the form (\ref{eqn_local_model}).
\end{enumerate}
The manifolds $E$ and $E'$, which are $\bbP^k$-bundles over $Y^s$, are called
the {\bf centers} of the transformation $\phi$ and its inverse, respectively.

\hfill

To summarize, Mukai's elementary transformation is the composition of a
contraction of a $\bbP^k$-bundle and of the inverse of another such contraction,
and transversally to the singular locus of the variety obtained by contraction
the transformation is modeled on the usual Mukai flop of the total space of the cotangent bundle of $\bbP^k$
in its zero section, as in \ref{defn_Mukai_flop}.

\subsection{Hu--Yau transformations and the Hu--Yau conjecture}\label{sec_HY}

In order to study bimeromorphic transformations of hyperk\"ahler manifolds
of arbitrary dimension, we will need to use bimeromorphic maps more general
than the Mukai's elementary transformations introduced in \ref{defn_Mukai_flop2}.
Following the idea of Hu and Yau \cite{_HuYau:AdvTeor_}, we will allow arbitrary bimeromorphic
transformations that generically look like Mukai flops in codimension two,
at the same time ignoring what happens in codimension greater than two.
We formalize this by introducing the following notion.

\hfill

\definition\label{defn_Hu_Yau}
Let $\phi\colon X\dra X'$ be a bimeromorphic map between two complex manifolds.
We will say that $\phi$ is a {\bf Hu--Yau transformation} if there exist closed
subvarieties $Z\subset X$ and $Z'\subset X'$ of codimension strictly greater than two,
such that over the open subsets $X_0 = X\setminus Z$ and $X_0'= X'\setminus Z'$
$$\phi|_{X_0}\colon X_0\dra X_0'$$ is a Mukai elementary transformation
in codimension two in the sense of \ref{defn_Mukai_flop2}.

\hfill

We reproduce the following conjecture made by Hu and Yau in \cite[Conjecture 7.3]{_HuYau:AdvTeor_},
stating it only for compact hyperk\"ahler manifolds, since this will be the only case
considered in the present paper. The definition of a compact hyperk\"ahler manifold
is recalled below, in Section \ref{sec_HK}.

\hfill

\conjecture 
Let $\phi:\; M\dra M'$ be a bimeromorphic map of compact
hyperk\"ahler manifolds. Then, after removing subvarieties
of codimension greater than two, $\phi$ is a sequence of
Mukai's elementary transformations.

\hfill

The above conjecture predicts a description of codimension two structure of bimeromorphic
maps between hyperk\"ahler manifolds. Since we are
allowed to remove subvarieties of codimension greater than two, the conjecture
does not predict anything about bimeromorphic maps that are regular in codimension two.
Our main result is the following statement. We use the notion of a Hu--Yau transformation
introduced above in \ref{defn_Hu_Yau}.

\hfill

\theorem\label{thm_main}
Let $M$ and $M'$ be compact hyperk\"ahler manifolds,
and $\phi:\; M \dra M'$ a bimeromorphic map that is not regular in codimension two. Then 
$\phi$ is a composition of a finite number of
Hu--Yau transformations between compact hyperk\"ahler manifolds.

\hfill

The proof of the above theorem will be given in Section \ref{sec_main_proof}.


\section{Hyperk\"ahler manifolds and MBM classes}\label{sec_HK}


In this section we fix the notation related to hyperk\"ahler manifolds that
will be used throughout the paper.

\subsection{Hyperk\"ahler manifolds and their Teichm\"uller spaces}

By a compact hyperk\"ahler manifold we will mean a compact simply connected complex manifold
that admits a hyperk\"ahler metric of maximal holonomy, see \cite{Hu}. To study bimeromorphic
maps between hyperk\"ahler manifolds it is sufficient to work with a fixed deformation class,
because by \cite[Proposition 27.8]{Hu} bimeromorphic hyperk\"ahler manifolds are deformation
equivalent. Therefore we fix the underlying $C^\infty$-manifold $M$. Henceforth, a hyperk\"ahler
manifold is the same thing as a complex structure of hyperk\"ahler type (i.e. admitting a hyperk\"ahler metric) on $M$.
Given such a complex structure $I$, we will write
$(M, I)$ for the corresponding hyperk\"ahler manifold. If the complex structure is clear from the
context or irrelevant, we will simply write $M$ instead of $(M, I)$.
The manifold $M$ is holomorphically symplectic and irreducible, i.e. the holomorphic symplectic form on $M$ is
unique up to multiplication by a constant. We will denote by $\sigma_I \in H^0(M, \Omega^2_M)$ a
symplectic form of unit volume, where the volume is defined as the integral of the form $\sigma_I^n\wdg\overline{\sigma}_I^n$.
Here $2n$ is the complex dimension of $M$. We will also assume that $b_2(M) \ge 6$.

We will denote by $q$ the Beauville-Bogomolov-Fujiki (BBF for short) quadratic form on $H^2(M, \bbQ)$.
Recall from \cite{Hu} that $q$ is of signature $(3, b_2(M) - 3)$
and satisfies the Fujiki relations:
$$
q(a)^n = c_M \int_M a^{2n}
$$
for any $a\in H^2(M, \bbQ)$, where $c_M\in \bbQ$ is a positive constant that can be chosen so
that $q$ is integral and primitive on $H^2(M, \bbZ)$, see \cite{So} for details.

Since $M$ is simply connected, its second homology is torsion-free,
hence $H_2(M, \bbZ)\simeq \mathrm{Hom}(H^2(M,\bbZ), \bbZ)$.
The form $q$ defines an embedding of lattices
$q\colon H^2(M, \bbZ) \hrarr H_2(M, \bbZ)$. This induces an embedding of $H_2(M, \bbZ)$ into the vector
space $H^2(X, \bbQ)$, in particular $q$ restricts to a rational quadratic form on the homology classes
of degree two. In what follows we will implicitly use this extension of $q$ to $H_2(M, \bbZ)$. 

The Teichm\"uller space of all complex structures of hyperk\"ahler type on $M$ will be denoted by $\Teich(M)$.
The point corresponding to the complex structure $I$ will be denoted by $[I]\in \Teich(M)$.
Recall the definition of the period domain:
$$
\DD = \{x\in \bbP H^2(M,\bbC) \st q(x) = 0, \, q(x,\bar{x}) > 0\}
$$
and the period map
$$
\rho\colon \Teich(M)\to \DD,
$$
sending a point $[I]\in \Teich(X)$ to $[\sigma_{I}]$. According to \cite{Hu}, if $(M, I_1)$ and
$(M, I_2)$ are bimeromorphic, then there exists a diffeomorphism $\alpha\colon M\to M$ such that the points $[I_1]$ and $[\alpha^*I_2]$
are non-separated in $\Teich(M)$.

\subsection{MBM classes and the K\"ahler cone}\label{sec_MBM}

Bimeromorphic hyperk\"ahler manifolds correspond to non-separated points in the Teichm\"uller
space, so for our purposes we may fix a connected component of the latter. This connected component
will be denoted by $\Teich^\circ(M)$. The non-separated points of the Teichm\"uller space
are distinguished by the K\"ahler cones of the corresponding manifolds. For $[I]\in \Teich(M)$
the K\"ahler cone of $(M, I)$, denoted by $\KK_I$, is the set of cohomology classes in
$H^{1,1}_I(M, \bbR)$ that can be represented by K\"ahler forms. To describe the structure
of the K\"ahler cone we will use the notion of an MBM class which we will now recall.

Given a cohomology class $\eta\in H^2(M, \bbQ)$, we define
$$
\Teich^\circ_\eta(M) = \{[I]\in \Teich^\circ(M)\st q(\rho(I), \eta) = 0\},
$$

in other words, the subspace of the complex structures where $\eta$ is of type $(1,1)$.

\hfill

\definition
A cohomology class $\eta\in H^2(M, \bbZ)$ is called {\bf an MBM class} if 
it satisfies the following conditions:
\begin{enumerate}
\item $\eta$ is primitive;
\item $q(\eta) < 0$;
\item for any $[I]\in \Teich^\circ_\eta(M)$ and any $\kappa\in \KK_I$ we have $q(\kappa, \eta) \neq 0$.
\end{enumerate}
The set of MBM classes on $M$ will be denoted by $\MBM(M)$ (see \cite{AV-MBM} for various equivalent characterizations of MBM classes).

\hfill

The notion of an MBM class has a deep connection with the geometry of rational curves on $M$.
We will partially recall this below, but see \cite[section 3]{_AV:Contraction_}
for a more detailed discussion. Recall that by \cite[Theorem 3.17]{_AV:Orbits_} the mapping class
group of $M$ acts on $\MBM(M)$ with finitely many orbits, in particular the BBF squares of the MBM
classes are bounded.

Given $[I]\in \Teich^\circ(M)$, recall that the restriction of $q$ to $H^{1,1}_I(M, \bbR)$
has signature $(1, b_2(M) - 3)$, therefore the set $\{x\in H^{1, 1}_I(X, \bbR)\st q(x) > 0\}$
has two connected components. We define the positive cone $\Pos_I$ to be the connected component
that contains the K\"ahler cone $\KK_I$. Given $[I']\in \Teich^\circ(M)$ and a
bimeromorphic map $\varphi\colon (M, I) \dra (M, I')$, recall from \cite{Hu} that
$\phi^*$ is a well defined map on $H^2(M, \bbC)$. The set $\varphi^*\KK_{I'}$
is an open convex cone contained $\CC_I$.

\hfill

\definition The union of the cones $\varphi^*\KK_{I'}$ for all bimeromorphic maps $\phi$ as above
is called {\bf the birational K\"ahler cone} of $(M, I)$, denoted by $\BK_I$.

\hfill

Note that the birational K\"ahler cone is not convex in general. To describe its structure we
introduce the wall and chamber decomposition of the positive cone. 
Given $\eta\in \MBM^{1,1}_I(X)$, denote by $W_\eta$ the corresponding wall in the positive cone,
i.e. the intersection of $\CC_I$ with the $q$-orthogonal complement of $\eta$. Define
$$
\CC_I^\circ = \CC_I\setminus \left(\cup_{\eta\in \MBM^{1,1}_I(X)} W_\eta\right).
$$

\hfill

\definition In the notation introduced above, {\bf a K\"ahler chamber} is a connected component of 
$\CC_I^\circ$.

\hfill

The following statement summarizes a description of the cones introduced above.
It is a combination of results proven in \cite{AV-MBM}, \cite{_AV:Orbits_}, \cite{_Markman:Survey_} and \cite{Bo}.

\hfill

\theorem\label{thm_chambers}
Assume that $b_2(M)\ge 5$ and let $[I]\in\Teich^\circ(M)$. Then we have the following:
\begin{enumerate}
	\item The collection of walls $W_\eta$, where $\eta\in \MBM^{1,1}_I(X)$, is locally finite in $\CC_I$\footnote{Remark that  
		walls may accumulate towards the the border of $CC_I$, the isotropic cone.} ;
\item The K\"ahler cone $\KK_I$ is a connected component of $\CC_I^\circ$;
\item For every K\"ahler chamber there exists a unique $[I']\in\Teich^\circ(M)$ such that
$\rho(I') =\rho(I)$ and $\KK_{I'}$ is equal to that chamber.
The manifolds $(M, I)$ and $(M, I')$ are bimeromorphic. This gives a bijection between the K\"ahler chambers
and the points $[I']\in\Teich(M)$ with $\rho(I) = \rho(I')$ and $(M, I')$ bimeromorphic to $(M, I)$;
\item The birational K\"ahler cone $\BK_I$ is the union of some of the K\"ahler chambers.
The closure $\overline{\BK_I}$ is the cone dual to the pseudo-effective cone of $M$, in particular
$\overline{\BK_I}$ is convex.
\end{enumerate}

We close this section by outlining the connection between the MBM classes and the structure
of extremal contractions of hyperk\"ahler manifolds. For more details we refer to \cite{_AV:Contraction_}.

\hfill

\definition Given a class $\eta\in \MBM^{1,1}_I(M)$ we will say that $\eta$
{\bf is extremal}, or that $\eta$ defines {\bf a wall of the K\"ahler cone}, if
$W_\eta\cap\overline{\KK_I}$ contains a non-empty open subset of $W_\eta$. 

\hfill

Any MBM class is extremal on some deformation of the ambient complex structure; in a given one,
MBM classes are extremal up to monodromy and birational transformation, see \cite{AV-MBM}.

By a rational curve $C$ in a hyperk\"ahler manifold $M$ we will mean the image
of a non-constant holomorphic map from $\bbP^1$ into $M$. Its homology class will
be denoted by $[C]$. In the following definition we implicitly use the embedding
$H^2(M,\bbZ)\hrarr H_2(M, \bbZ)$ given by the BBF form, as was explained above.

\hfill

\definition Let $\eta\in \MBM^{1,1}_I(M)$ be an extremal MBM class.
The corresponding {\bf full MBM locus} $E_\eta$ is the closure of the union
of all rational curves $C$ such that $[C]$ is proportional to $\eta$.

\hfill

Note that the full MBM locus need not be irreducible (or even connected).
We will need to talk about the codimension of the MBM loci. By this we
will mean the following.

\hfill

\definition Let $\eta\in \MBM^{1,1}_I(M)$ be an extremal MBM class. The {\bf codimension} of $\eta$
(and of $E_\eta$) is the minimum of the codimensions of the irreducible components
of $E_\eta$.  We will say that $\eta$ is {\bf divisorial} if it is of codimension one,
otherwise $\eta$ will be called {\bf non-divisorial}.

\hfill

\lemma Let $\eta$ be a divisorial MBM class. Then $E_\eta$ is an irreducible divisor.

\hfill

\begin{proof}
By the definition of the codimension of $\eta$, $E_\eta$ contains at least one irreducible component
of codimension one, call this component $D$. The divisor $D$ is uniruled, covered by
deformations of a rational curve $C$ with $[C]$ proportional to $\eta$, by the definition
of $E_\eta$. It is known (see e.g. \cite[Lemma 2.11]{_AV:MBM_}) that $[C]$ is proportional
to $[D]$, therefore $[D]$ is proportional to $\eta$. Assume that $E'$ is another irreducible component of $E_\eta$.
Then $E'$ contains a rational curve $C'$ that does not lie inside of $D$
and with $[C']$ proportional to $\eta$. Since $C'$ is not contained in $D$, their intersection
index is non-negative, contradicting the fact that $q(\eta) < 0$.   
\end{proof}

\hfill

In the terminology of \cite{Bo}, the divisorial MBM loci are the prime exceptional divisors, and
their cohomology classes span the extremal rays of the pseudo-effective cone of $M$. 
Therefore the divisorial MBM classes define the walls of $\overline{\BK_I}$, i.e.
for a divisorial MBM class $\eta$ the intersection $W_\eta\cap \overline{\BK_I}$
contains a non-empty open subset of $W_\eta$ that lies on the boundary of $\overline{\BK_I}$.
The non-divisorial extremal MBM classes define the internal walls that cut the birational
K\"ahler cone into K\"ahler chambers.

We close the section by recalling the following result (see also \cite[Corollary 5.9]{BL}).

\hfill

\theorem\cite[Theorem 4.6]{_AV:Contraction_}\label{thm_AV} Let $\eta\in \MBM^{1,1}_I(M)$ be an extremal MBM class
on a hyperk\"ahler manifold with $b_2(M)\ge 5$. Then there exists a normal symplectic variety $N$
and a proper bimeromorphic morphism $\pi\colon M\to N$ such that the exceptional
set of $\pi$ is the full MBM locus $E_\eta$.


\section{Wall-crossing flops}\label{sec_WC}


In this section we study the bimeromorphic transformations corresponding to
passing through a wall $W_\eta$ defined by an extremal non-divisorial MBM class $\eta$.

Consider a point $[I]\in \Teich^\circ(M)$ and the corresponding birational K\"ahler cone $\BK_I$.

\hfill

\definition\label{defn_wc_flop}
We will say that two chambers $\KK_1, \KK_2\subset \BK_I$ are {\bf adjacent} if 
the intersection $\overline{\KK_1} \cap \overline{\KK_2}$ contains a non-empty open
subset of one of the walls $W_\eta$ for some $\eta\in \MBM^{1,1}_I(M)$.

\hfill

\definition
Let $\varphi\colon M_1 \dra M_2$ be a bimeromorphic map between hyperk\"ahler manifolds
with K\"ahler cones $\KK_1$ and $\KK_2$. If $\KK_1$ and $\varphi^*\KK_2$ are adjacent,
then $\varphi$ is called a {\bf wall-crossing flop}.

\hfill

It will be convenient for us to work with wall-crossing flops that satisfy
an additional condition, namely that $\phi^*$ acts as the identity on the second
cohomology of $M$. It is always possible to arrange this due to the following statement.

\hfill

\proposition\label{prop_chambers}
Assume that $\KK'\subset \BK_I$ is a chamber.
Then there exists a complex structure $I'$ of hyperk\"ahler type on $M$, unique up to isotopy, and a bimeromorphic
map $\varphi\colon (M, I) \dra (M, I')$ such that $\varphi^*$ is the identity map on $H^2(M, \bbZ)$
and $\KK' = \KK_{I'}$.

\hfill

\begin{proof} By the definition of $\BK_I$, there exists $[\tilde{I}]\in \Teich^\circ(M)$ and
a bimeromorphic map $\psi\colon (M, I)\dra (M, \tilde{I})$ such that $\KK' = \psi^*\KK_{\tilde{I}}$.
By \cite[Proposition 27.8]{Hu} there exist two families $\pi_1\colon \MM\to \Delta$ and $\pi_2\colon \tilde{\MM}\to \Delta$
of hyperk\"ahler manifolds over the unit disc and a bimeromorphic map $\Psi\colon \MM\dra \tilde{\MM}$
with the following properties. First, we have $\MM_0 \simeq (M, I)$ and $\tilde{\MM}_0\simeq (M,\tilde{I})$, where
$\MM_t$ and $\tilde{\MM}_t$ denote the fibers over $t\in\Delta$. Second, $\Psi$ induces a fiberwise isomorphism
over $\Delta\setminus\{0\}$, and $\Psi_0\colon \MM_0\dra \tilde{\MM}_0$ equals $\psi$.

Trivializing the families $\MM$ and $\tilde{\MM}$ by Ehresmann's theorem, we may identify $\MM_{t_0}$ and $\tilde{\MM}_{t_0}$
with $M$ as $C^\infty$-manifolds, where $0\neq t_0\in \Delta$ is some base point. Then $\Psi_{t_0}$ defines a diffeomorphism
of $M$. By the second property mentioned above, the action of $\Psi_{t_0}$ on $H^2(M, \bbZ)$ coincides with the action of $\psi$.

Define $I' = \Psi_{t_0}^*\tilde{I}$. Then $\Psi_{t_0}^{-1}\colon (M, \tilde{I}) \to (M, I')$ is a biholomorphic isomorphism,
hence $\KK_{\tilde{I}} = (\Psi_{t_0}^{-1})^*\KK_{I'}$. Define $\phi = \Psi_{t_0}^{-1}\circ \psi$. Then by construction
$\phi^*$ acts on $H^2$ as the identity. The claim about uniqueness of the complex structure $I'$ up to isotopy follows
from the description of the fibers of the period map, as was recalled in \ref{thm_chambers}.
\end{proof}

\hfill

\corollary\label{cor_flop_identity}
Let $\phi\colon M_1\dra M_2$ be a bimeromorphic map, where $M_1 = (M, I_1)$ and $M_2 = (M, I_2)$.
Then there exists a unique $[I_2']\in \Teich^\circ(M)$ and an isomorphism $\xi\colon M_2 \to M_2'$,
where $M_2' = (M, I_2')$, such that $\phi' = \xi\circ \phi\colon M_1\dra M_2'$ is a bimeromorphic
map that acts as the identity on the second cohomology of $M$. If $\phi$ is a wall-crossing flop,
then so is $\phi'$.

\hfill

\begin{proof}
Let $\KK' = \phi^*\KK_{I_2}$. By \ref{prop_chambers} there exists a complex structure $I_2'$,
unique up to isotopy, and a bimeromorphic map $\varphi'\colon M_1 \dra M_2' = (M, I_2')$
such that $\varphi'$ acts as the identity on $H^2$ and and $\KK' = \KK_{I_2'}$. The
composition $\phi\circ (\phi')^{-1}$ is a bimeromorphic map from $M_2'$ to $M_2$, and
it maps the K\"ahler cone of $M_2$ to the K\"ahler cone of $M_2'$, hence by \cite[Proposition 27.6]{Hu}
it is biholomorphic. Defining $\xi$ to be the inverse of this map proves the first claim.
If $\phi$ is a wall-crossing flop, then by \ref{defn_wc_flop} $\KK'$ is adjacent to
the K\"ahler cone of $M_1$. By construction $(\phi')^*\KK_{I'_2} = \KK'$, hence $\phi'$
is also a wall-crossing flop.
\end{proof}

\hfill

\lemma\label{cor_flop_difference}
Let $\phi\colon M_1 \dra M_2$ and $\phi'\colon M_1 \dra M_2$ be two bimeromorphic
maps acting as the identity on the second cohomology. Then $\phi' = \xi\circ \phi$,
where $\xi$ is a biholomorphic automorphism of $M_2$ that acts as the identity on the second cohomology.

\hfill

\begin{proof}
Define $\xi = \phi'\circ\phi^{-1}$ and conclude that it is a biholomorphic map by \cite[Proposition 27.6]{Hu}.
\end{proof}

\hfill

Now can describe the nature of wall-crossing flops, showing that they are compositions of
a blow-down followed by a blow-up, similarly to Mukai's elementary transformations.

\hfill

\proposition\label{prop_flop}
Let $\varphi\colon M_1 \dra M_2$ be a wall-crossing flop. If $M_i$ are non-projective,
assume additionally that $b_2(M_i) \ge 6$. Then there exists a compact
symplectic variety $N$ and two small contractions $\pi_1\colon M_1 \to N$ and $\pi_2\colon M_2\to N$
such that $\varphi = \pi_2^{-1}\circ \pi_1$.

\hfill

The proof of the above proposition will be divided in two parts, depending on whether $M_i$ are projective.

\subsection{The projective case}

We prove \ref{prop_flop} assuming that $M_i$ are projective. Let $M_1 = (M, I_1)$ and $M_2 = (M, I_2)$.
By \ref{cor_flop_identity} we may assume that the map $\varphi$ acts as the identity on $H^2(M, \bbZ)$.
By the definition of a wall-crossing flop, there exists a wall $W_z$ of the positive cone $\CC_{I_1}$ such
that the intersection $\overline{\KK_{I_1}} \cap \overline{\KK_{I_2}}$ contains a non-empty open subset of $W_z$. Since the
manifolds $M_i$ are projective, this open subset contains an integral element $c\in H^2(M, \bbZ)$. The class
$c$ is of Hodge type $(1, 1)$ both on $M_1$ and $M_2$, so there exist $L_1\in\Pic(M_1)$ and $L_2\in\Pic(M_2)$
with $c_1(L_i) = c$. By construction, the class $c$ is big and nef both on $M_1$ and $M_2$, so by
a theorem of Kawamata \cite[Theorem 6.1]{Kw} the bundles $L_i$ are semiample
and define small contractions $f_1\colon M_1\to N_1$ and $f_2\colon M_2\to N_2$ onto projective symplectic
varieties $N_1$ and $N_2$. Let $\psi\colon N_1 \dra N_2$ be the composition $f_2\circ \varphi\circ f_1^{-1}$.
Note that $L_i \simeq f_i^*\OO_{N_i}(1)$, where $\OO_{N_i}(1)$ are very ample line bundles on $N_i$, and
$\varphi^* L_2 \simeq L_1$ by construction. We see that $\psi^*\OO_{N_2}(1) \simeq \OO_{N_1}(1)$,
therefore by the Matsusaka--Mumford criterion \cite[Theorem 11.39]{Ko} we conclude that $\psi$ is an
isomorphism. Identifying $N_1$ with $N_2$ via $\psi$ and letting $N=N_1=N_2$, we finish the proof.

\subsection{The non-projective case}

To prove \ref{prop_flop} in the case when $M_i$ are non-projective we will use ergodic
theory, as in \cite{_AV:Contraction_}. We first make the following observation.

\hfill

\lemma\label{lem_abs_auto}
Let $\pi\colon M\to N$ be an extremal MBM contraction with the exceptional set $E_\eta$,
as in \ref{thm_AV}, and $\xi\colon M\to M$ an automorphism that acts trivially on $H^2$.
Then $\xi$ descends to an automorphism of $N$.

\hfill

\begin{proof}
Consider the contraction $\pi' = \pi\circ \xi$. Since $\xi$ acts trivially on the second cohomology,
it preserves homology classes of the rational curves that cover the full MBM locus $E_\eta$.
By the definition of a full MBM locus, $\pi'$ preserves $E_\eta$ and contracts all the fibers
of $\pi$. Consider the map $\pi\times\pi' \colon M \to N\times N$ and let $Z$ be its
image. The variety $Z$ is finite of degree one over both factors of $N\times N$, and since
$N$ is normal, $Z$ is isomorphic to $N$, therefore being the graph of an automorphism of $N$.
\end{proof}

\hfill

Fix $\eta\in \MBM^{1,1}_I(M)$ and recall that given a K\"ahler chamber $\KK\subset \CC_I^\circ$
we say that $\eta$ defines a wall of $\KK$ if $\overline{\KK}\cap W_\eta$ contains a non-empty
open subset of the wall $W_\eta$. Define the following open subspaces of $\Teich_\eta^\circ(M)$:
\begin{eqnarray}
\Teich^+_\eta(M) = \{ [I]\in \Teich^\circ_\eta(M)\st \eta \mbox{ defines a wall of } \KK_I \mbox{ and } q(\eta, \KK_I) > 0\}\nonumber\\
\Teich^-_\eta(M) = \{ [I]\in \Teich^\circ_\eta(M)\st \eta \mbox{ defines a wall of } \KK_I \mbox{ and } q(\eta, \KK_I) < 0\}\nonumber
\end{eqnarray}

Note that there is a one-to-one correspondence between the points of the two spaces above: given $[I^+]\in \Teich^+_\eta(M)$,
the intersection $\overline{\KK}_{I^+} \cap W_\eta$ contains a non-empty open subset of $W_\eta$, so there exists a unique
chamber ${\KK}_{I^-}$ on the opposite side of the wall $W_\eta$, so that $\overline{\KK}_{I^-} \cap W_\eta = \overline{\KK}_{I^+} \cap W_\eta$
and $[I^-]\in \Teich^-_\eta(M)$. The two K\"ahler cones ${\KK}_{I^\pm}$ are adjacent and by \ref{prop_chambers} there exists
a wall-crossing flop $\phi\colon (M, I^+)\dra (M, I^-)$. It is clear that all wall-crossing flops across the walls defined by $\eta$ arise this way.

We have two period maps $\rho^\pm\colon \Teich^\pm_\eta(M)\to \DD_\eta$, where $\DD_\eta = \DD\cap \bbP(\eta^\perp)$.
We see that the fibers of $\rho^\pm$ over a point $[\sigma_I]\in \DD_\eta$ are in one-to-one correspondence
with the chambers of the wall $W_\eta = \CC_I \cap \eta^\perp$. Those chambers are cut out by the orthogonal complements
to the MBM classes in $\MBM^{1,1}_I(M)$ different from $\pm \eta$.

Denote by $G$ the mapping class group of $M$ and by $G_\eta\subset G$ the stabilizer of $\eta$. Recall from \cite[Theorem 2.14]{_AV:Contraction_}
the classification of $G_\eta$-orbits in $\DD_\eta$. It follows from the classification that given a period $[\sigma_I]\in \DD_I$
of a non-projective hyperk\"ahler manifold $(M,I)$, its orbit $G_\eta\cdot[\sigma_I]$ contains a period of some projective
hyperk\"ahler manifold in its closure.

Let us finish the proof of \ref{prop_flop} in the non-projective case. Let $\phi\colon M^+\dra M^-$
be a wall-crossing flop that acts trivially on $H^2$, where $M^\pm = (M, I^\pm)$ for some $[I^\pm]\in \Teich^\pm_\eta(M)$.
Let $\pi^\pm\colon M^\pm \to N^\pm$ be the extremal contractions that exist by \ref{thm_AV}.
We need to show that the bimeromorphic map $N^+\dra N^-$ induced by $\phi$ is an isomorphism.
Assuming that $I^\pm$ are non-projective, we use the $G_\eta$-orbit classification \cite[Theorem 2.14]{_AV:Contraction_}
and \cite[Theorem 3.1]{V2} to find a sequence of diffeomorphisms
$\gamma_i$ such that $[\gamma_i^*I^+] \to [I_0^-]$, $i\to \infty$, for some projective 
complex structures $[I^\pm_0]\in \Teich^\pm_\eta(M)$. Clearly, also $[\gamma_i^*I^-]\to [I_0^-]$. Denote $M_0^\pm = (M, I_0^\pm)$ and let
$\phi_0\colon M_0^+ \dra M_0^-$ be a wall-crossing flop acting trivially on $H^2$. By the projective
case of the proposition proven above, we can decompose $\phi_0$ into the composition of a contraction
$\pi^+_0\colon M_0^+ \to N_0$ and the inverse of a contraction $\pi^-_0\colon M_0^-\to N_0$.

We refer to \cite{BL} and references therein for the notion of a locally trivial deformation of a complex variety.
Let $\nu\colon \NN\to B$ be the universal locally trivial deformation of $N_0$ over a polydisc, and $\Pi^\pm\colon \MM^\pm\to \NN$
the corresponding deformations of the contractions $\pi_0^\pm$ that exist by \cite[Theorem 1.1]{BL}.
The morphisms $\Pi^\pm$ induce a family $\Phi\colon \MM^+ \dra \MM^-$, $\Phi_0=\phi_0$,
of wall-crossing flops (each $\Phi_t$ over $\NN_t$), where $\mu^\pm\colon \MM^\pm\to B$ are the
deformations of $M_0^\pm$ preserving the contractions. The latter induce embeddings $B\hrarr \Teich^\pm_\eta(M)$, making $B$ an open
neighborhood of $[I_0^\pm]$ in $\Teich^\pm_\eta(M)$. The convergence $[\gamma_i^*I^\pm] \to [I_0^\pm]$, $i\to \infty$
implies that we can find $t\in B$ and $\gamma\in G_\eta$ such that $\MM_t^\pm \simeq (M, \gamma^*I^\pm)$.

Now we have two wall-crossing flops $\Phi_t$ and  $\gamma^{-1}\circ \phi\circ\gamma \colon \MM_t^+ \dra \MM_t^-$.
By \ref{cor_flop_difference} we have $\gamma^{-1}\circ \phi\circ\gamma = \xi\circ\Phi_t$ for some automorphism $\xi$
of $\MM_t^-$ acting trivially on $H^2$.
By construction, the flop $\Phi_t$ induces an isomorphism between the bases $\NN_t^\pm$ of the extremal contractions
$\Pi^\pm_t\colon \MM_t^\pm\to \NN_t^\pm$. By \ref{lem_abs_auto} the map $\xi$ induces an automorphism of $\NN_t^-$.
Therefore the flop $\gamma^{-1}\circ \phi\circ\gamma$ induces an isomorphism between $\NN_t^\pm$,
and the same is true about $\phi$, since $\gamma\colon \MM_t^\pm\to (M, I^\pm)$ is an isomorphism.
This completes the proof of \ref{prop_flop}.


\subsection{Decomposition}


We now show how to decompose an arbitrary bimeromorphic map into a composition of wall-crossing flops.

\hfill

\proposition\label{prop_decomp}
Assume that $\varphi\colon M \dra M'$ is a bimeromorphic map between compact hyperk\"ahler
manifolds with $b_2 \ge 5$. Then there exists a finite sequence of wall-crossing flops $\varphi_i\colon M_i\dra M_{i+1}$,
$i=0,\ldots,m-1$ such that $M_0 = M$, $M_m = M'$ and $\varphi = \varphi_m \circ \ldots \circ \varphi_0$,

\hfill

\begin{proof}
By \ref{prop_chambers} we may assume that $\varphi\colon (M, I)\dra (M, I')$ acts as the identity
on $H^2$ and $\KK_I$, $\KK_{I'}$ are two chambers of the birational K\"ahler cone $\BK_I$. Choose
two points $a\in \KK_I$, $b \in \KK_{I'}$ and connect them by a straight segment $\ell\subset H^{1,1}_I(M, \bbR)$.
By \cite[Theorem 3.17]{_AV:Orbits_} and \cite[Proposition 5.14]{SSV} the collection of walls $W_z$, $z\in \MBM^{1,1}_I(M)$
is locally finite in $\CC_I$, and since the segment $\ell$ is compact, it intersects only a finite number of
walls $W_{z_1},\ldots,W_{z_m}$. Moreover, since the chambers $\KK_{I}$ and $\KK_{I'}$ are open, we may,
by perturbing the endpoints $a$ and $b$, make sure that all points of intersection $W_i\cap \ell$ are distinct.
In this case, whenever the segment $\ell$ crosses a wall $W_{z_i}$, it passes to an adjacent chamber.
We therefore get a sequence of adjacent chambers $\KK_0 = \KK_I, \KK_1,\ldots,\KK_m = \KK_{I'}$. For every
chamber $\KK_j$ there exists a complex structure $I_j$ such that $\KK_{I_j} = \KK_j$ and a wall-crossing
flop $\varphi_j\colon (M, I_j)\dra (M,I_{j+1})$ acting as the identity on $H^2$. Consider the composition
$\alpha = \varphi^{-1} \circ \varphi_m \circ \ldots \circ \varphi_0$, so that $\alpha\colon (M, I)\dra (M, I)$.
By construction, the map $\alpha$ acts as the identity on $H^2$, hence it preserves the K\"ahler cone $\KK_I$,
therefore by \cite[Proposition 27.6]{Hu} the map $\alpha$ is a regular automorphism. Replacing $\varphi_0$ with $\varphi_0\circ \alpha^{-1}$,
we finish the proof.
\end{proof}


\subsection{Wall-crossing flops in codimension two}


We apply the main theorem of \cite{WW} to establish the structure of a wall-crossing flop
in codimension two. The main step is to reduce the question to the case of four-dimensional symplectic
manifolds. We do this by taking a slice transversal to the indeterminacy locus of the flop.

Assume that $M$ is a holomorphically symplectic manifold and $\pi\colon M\to N$
a proper bimeromorphic map onto a normal variety $N$. We recall that according to \cite{Kl} the variety
$N$ is symplectic and admits a unique stratification $N = \coprod N_i$, where $N_i$ are locally closed
non-singular symplectic subvarieties of $N$. We will call it {\bf the Kaledin stratification}.
Recall also that the symplectic forms on $M$ and $N$ induce Poisson brackets on their structure
sheaves, so that $\pi$ is a morphism of Poisson varieties.

The following proposition essentially follows from \cite[Theorem 2.3]{Kl}. In loc. cit., both the slice of $N$
transversal to a stratum of the Kaledin stra\-tification and a local product decomposition are constructed on the level of formal germs.
Our proposition contains two minor modifications: first, our slices and the decomposition are constructed locally in the analytic topology;
second, we construct them in the presence of a crepant resolution, and show that the resolution also
admits a local product decomposition. The proof is the same as in loc. cit., we give it here for the sake of completeness.

\hfill

\proposition\label{prop_slice}
Let $M$ be a holomorphically symplectic manifold
and $\pi\colon M \to N$ a proper small contraction onto a normal symplectic variety.
Let $Z\subset N$ be a stratum of the Kaledin stratification.
Then for any point $z\in Z$ there exists an open subset $U\subset N$ containing $z$
and a closed symplectic subvariety $N'\subset U$ called {\bf the Kaledin slice} such that:
\begin{enumerate}
\item $z\in N'$ and $\dim(N') = \mathrm{codim}_N(Z)$;
\item $M' = \pi^{-1}(N')$ is a non-singular symplectic variety and $\pi$ induces a small contraction of $M'$ onto $N'$;
\item The subset $U$ and its preimage in $M$ admit product decompositions as symplectic varieties:
$$
U \simeq N' \times \Delta^{2d}, \quad \pi^{-1}(U) \simeq M' \times \Delta^{2d},
$$
where $2d$ is the dimension of $Z$ and the polydisc $\Delta^{2d}$ is endowed with the standard symplectic structure.
\end{enumerate}

\hfill

\begin{proof}
Let $E\subset M$ be the exceptional locus of $\pi$. Since $\pi$ is a small con\-traction,
all irreducible components of $E$ are of codimension at least two in $M$.
We will prove the proposition by induction on $d$. The case $d = 0$ being trivial, we assume that $d \ge 1$.

We let $U$ be a small analytic neighborhood of $z$ in $N$ that we may shrink when necessary.
Given a holomorphic function $f\in \OO(U)$, we will denote by $H_f$ its Hamiltonian vector field,
i.e. the vector field corresponding to the derivation $g\mapsto \{f,g\}$, where $\{\cdot,\cdot\}$ is 
the Poisson structure on $U$ coming from the symplectic form. We will also let $\tilde{f}$ be the
pullback of $f$ to $\pi^{-1}(U)$ and $H_{\tilde{f}}$ its Hamiltonian vector field. Observe
that since $\pi$ is a morphism of Poisson varieties, we have $d\pi(H_{\tilde{f}}) = H_{f}$.

Given a tangent vector $v\in T_z Z$, we claim that there exists a holomorphic function $f\in \OO(U)$
such that $H_{f}(z) = v$. Indeed, $Z$ is a smooth manifold with a symplectic form $\sigma$ induced
by the Poisson structure on $N$. We can therefore find a function $g\in \OO(Z\cap U)$ such that
$dg_z = v\cnv\sigma$. The function $g$ is a restriction of some $f\in \OO(U)$, and since $Z$ is
a Poisson submanifold of $N$, we have by construction $H_f(z) = v$.

Given a function $f$ as above, denote by $N_f$ its zero locus, and by $M_f$ the preimage of $N_f$
in $M$, i.e. the zero locus of $\tilde{f} = f\circ\pi$. Observe that after shrinking $U$ we may
assume that $H_f$ is nowhere vanishing, and hence the same is true of $H_{\tilde{f}}$. It follows that
$M_f$ is a non-singular subvariety of $\pi^{-1}(U)$, indeed, $df$ and $d\tilde{f}$ are non-zero by duality. Given a point $x\in E\cap \pi^{-1}(U)$, we have
$d\pi(H_{\tilde{f}}(x)) = H_f(\pi(x))$ and since $H_f$ is tangent to the strata of the Kaledin stratification,
the vector $H_{\tilde{f}}(x)$ is tangent to $E$. Moreover, $H_{\tilde{f}}(x)$ is not contained in
the kernel of $d\pi$, therefore $E_f = E\cap M_f$ is of codimension one in $E$.

We now choose two linearly independent vectors $v_1, v_2\in T_z Z$ with $\sigma(v_1, v_2)$ = 1,
and the corresponding functions $f_1$ and $f_2$ as above. Shrinking $U$ we may assume that the function
$\{f_1, f_2\}$ is nowhere vanishing in $U$, so that $H_{f_1}$ and $H_{f_2}$ are linearly
independent in $U$. Then the hypersurfaces $M_{f_1}$ and $M_{f_2}$ intersect transversally
and their intersection $M''$ is a submanifold of $\pi^{-1}(U)$. The tangent bundle of $M''$ is
the symplectic orthogonal complement to the span of $H_{f_1}|_{M''}$ and $H_{f_2}|_{M''}$,
and since $\{f_1, f_2\}$ is nowhere vanishing, $M''$ is a symplectic submanifold of $\pi^{-1}(U)$.
Observe that by the analogous argument $E'' = E_{f_1} \cap E_{f_2}$ is of codimension two in $E\cap \pi^{-1}(U)$,
and $\pi|_{M''}$ is a small contraction onto $N'' = N_{f_1}\cap N_{f_2}$.

We may apply the induction hypothesis and obtain the product decompositions $N'' = N' \times \Delta^{2d - 2}$
and $M'' = M' \times \Delta^{2d - 2}$. To complete the proof, observe that
$\tilde{F} = (\tilde{f}_1, \tilde{f}_2)\colon \pi^{-1}(U) \to \bbC^2$ defines a submersion onto a neighborhood of the origin in $\bbC^2$.
The fiber over the origin is $M''$, and the submersion factors through $F = (f_1, f_2)\colon U \to \bbC^2$.
We may replace $U$ by the preimage of a small bidisc in $\bbC^2$ centered at the origin.
The vector fields $H_{f_1}$ and $H_{f_2}$ give a splitting of the exact sequence defining the relative tangent bundle
of the submersion $\tilde{F}$. So, we have $\pi^{-1}(U) \simeq M'' \times \Delta^2$ as soon as the vector
fields $H_{\tilde{f}_1}$ and $H_{\tilde{f}_2}$ commute, that is, generate a foliation. To achieve this we
need to modify the function $f_2$ without changing its zero locus.

We will replace $f_2$ with another function $f_2'$ such that $\{f_1, f_2'\} = 1$ in $U$
and $N_{f_2} = N_{f_2'}$. We do it as follows. On $\pi^{-1}(U)$, the vector field $H_{\tilde{f}_1}$ is
transversal to the hypersurface $M_{f_2}$, since by construction $\{f_1, f_2\}$ is nowhere
vanishing. Therefore there exists a unique function $\tilde{f}_2'$ in a neighborhood of $M_{f_2}$
such that $\tilde{f}_2'$ vanishes on $M_{f_2}$ and $\{\tilde{f}_1, \tilde{f}_2'\} = \LL_{H_{\tilde{f}_1}}(\tilde{f}_2') = 1$
(one can see this by using the local canonical form of a vector field in a neighborhood of its
regular point). After possibly shrinking $U$, the function $\tilde{f}_2'$ is the pullback of some $f_2'\in \OO(U)$, because
$\pi$ is an isomorphism in codimension one and $N$ is a normal variety.
Now we have $[H_{f_1}, H_{f_2'}] = H_{\{f_1,f_2'\}} = 0$, so the vector fields $H_{\tilde{f}_1}$ and $H_{\tilde{f}_2'}$
commute and define a regular foliation transversal to $M''$. The leaves of this foliation give the necessary splitting.
\end{proof}

\hfill

We next describe the structure of wall-crossing flops in codimension two,
showing that they are Mukai's elementary transformations in the sense of \ref{defn_Mukai_flop2}.

\hfill

\proposition\label{prop_WW}
Let $\phi\colon M_1 \dra M_2$ be a wall-crossing flop between hyperk\"ahler manifolds.
If $M_i$ are non-projective, assume additionally that $b_2(M_i) \ge 6$.
Then there exist two closed subvarieties $Z_1\subset M_1$ and $Z_2\subset M_2$ of codimension
strictly greater than two, such that $\phi\colon M_1\setminus Z_1\dra M_2\setminus Z_2$ is
Mukai's elementary transformation in codimension two (or an isomorphism).

\hfill

\begin{proof}
We deduce from \ref{prop_flop} that $\phi = (\pi_2)^{-1}\circ \pi_1$ where $\pi_i\colon M_i\to N$
are small contractions onto a symplectic variety $N$. Let $E_i\subset M_i$ be the exceptional loci
of the contractions. Since the MBM class defining the flop $\phi$ is non-divisorial, all irreducible
components of $E_i$ are of codimension strictly greater than one, and the singular locus of $N$ is of
codimension strictly greater than three. We are interested in the structure of the flop over the
codimension four strata of the Kaledin stratification.

Let $Z$ be the union of the strata of codimension six and higher (recall that the strata
of the Kaledin stratification are even-dimensional, being symplectic manifolds). Let 
$Z_i = \pi_i^{-1}(Z)$. The contractions $\pi_i$ being semismall, the subvarieties $Z_i$ are
of codimension at least three in $M_i$. Let $N^\circ = N\setminus Z$ and $M_i^\circ = M_i\setminus Z_i$.

Note that the singular locus of $N^\circ$ is either empty or the union of strata of codimension four in the Kaledin
stratification. If it is empty, by semismallness our flop is an isomorphism in codimension two. 
	Otherwise, given a singular point $z\in N^\circ$ we apply \ref{prop_slice} to construct the
Kaledin slice $N'$ passing through the point $z$. The slice $N'$ is four-dimensional and $z$ is
its isolated singular point. It has two crepant resolutions $M_1' = \pi_1^{-1}(N')$ and $M_2' = \pi_2^{-1}(N')$.
We now have to show that the preimages of $z$ in $M_i'$ are isomorphic to $\bbP^2$, since in
that case the product decomposition from \ref{prop_slice} implies that $\phi|_{M_1^\circ}$ is Mukai's elementary
transformation in codimension two.

We apply \cite[Theorem 1.1]{WW} to the small contractions $\pi_i\colon M_i' \to N'$ to conclude the proof.
Note that in loc. cit. the result is stated for quasi-projective varieties (however, see the remark following
the statement of Theorem 1.1. in loc. cit.). To reduce to the quasi-projective case, we use deformation invariance
of the contraction centers \cite[Theorem 1.8]{_AV:Contraction_} and deform $M_i$ preserving the wall that
defines $\phi$ to make $M_i$ projective. Next we note that in the proof of \ref{prop_slice} the stratum of the
Kaledin stratification is a quasi-projective manifold. In the induction step of the proof, to construct
the slice $N'$, we may choose the functions $f_1$ and $f_2$ to be algebraic, defined on a Zariski-open subset
of $N$. The slice $N'$ and its preimages $M_i'$ constructed this way will also be quasi-projective,
and we may apply \cite[Theorem 1.1]{WW} to them, concluding the proof.
\end{proof}

\hfill

The following statement follows directly from the above proposition and \ref{defn_Hu_Yau}.

\hfill

\corollary\label{cor_HU}
Any wall-crossing flop between hyperk\"ahler manifolds (assuming $b_2 \ge 6$ in the non-projective case)
that is not regular in codimension two is a Hu--Yau transformation.


\section{The symplectic rank and bimeromorphic maps}


\subsection{Definition of the symplectic rank}

Let $M$ be a holomorphic symplectic manifold with symplectic form $\sigma$
and $Z\subset M$ an irreducible subvariety. If $z\in Z$ is a smooth point,
then the rank of the two-form $\sigma|_{T_{Z,z}}$ is an even integer. This integer
is constant on a dense Zariski-open subset of $Z$, i.e. for a general point of $Z$. 
\begin{defn}
For an irreducible subvariety $Z\subset M$ its symplectic rank
$\rk_\sigma(Z)$ is the rank of the form $\sigma|_{T_{Z,z}}$ for a general smooth point $z\in Z$.
For a subvariety $Z =\cup_i Z_i$ with irreducible components $Z_i\subset M$ define
$\rk_\sigma(Z) = \max_i{\rk_\sigma(Z_i)}$.
\end{defn}

We record a few useful properties of the symplectic rank.

\hfill

\proposition\label{prop_srank}
Let $M$ be a holomorphic symplectic manifold with symplectic form $\sigma$
and $Z\subset M$ an irreducible subvariety.
\begin{enumerate}
\item $Z$ is isotropic if and only if $\rk_\sigma(Z) = 0$;
\item We always have $\rk_\sigma(Z) \ge \dim(Z) - \mathrm{codim}(Z)$
and the equality holds if and only if $Z$ is coisotropic;
\item If $Z'\subset Z$ is another irreducible subvariety, then $\rk_\sigma(Z')\le \rk_\sigma(Z)$;
\item If $W$ is an irreducible variety and $f\colon W\to Z$ a surjective
morphism, then $\rk_\sigma(Z)$ is the rank of $f^*\sigma$ at a general smooth point of $W$.
\end{enumerate}

\hfill

\begin{proof}
{\it (1), (2) and (3):} Obvious from the analogous statements about a subspace in a symplectic vector space.

{\it (4):} Resolving the singularities of $W$ we may assume that $W$ is smooth. Then the morphism
$f$ is generically smooth over the smooth locus of $Z$. We may find a point $w\in W$ such that
$df_w\colon T_{W,w}\to T_{Z, f(w)}$ is surjective and $\sigma|_{T_{Z, f(w)}}$ is of rank $\rk_\sigma(Z)$.
Then $f^*\sigma|_{T_{W,w}}$ is of the same rank as $\sigma|_{T_{Z, f(w)}}$. The points $w$ with this
property form a non-empty open subset of $W$, hence the claim.
\end{proof}

\hfill

We next observe that the symplectic rank is invariant under bimeromorphic transformations.

\hfill

\proposition\label{prop_srank_inv}
Let $\varphi\colon M\dra M'$ be a bimeromorphic map between hyperk\"ahler manifolds
with symplectic forms $\sigma$ and $\sigma'$,
$Z\subset M$ a closed subvariety and $Z' = \varphi(Z)$ its full transform under $\varphi$.
Then $\rk_\sigma(Z) = \rk_{\sigma'}(Z')$.

\hfill

\begin{proof}
Resolve the indeterminacies of $\phi$ by a bimeromorphic morphism $\pi_1\colon \tilde{M}\to M$,
so that $\phi = \pi_2\circ\pi_1^{-1}$ for some bimeromorphic morphism $\pi_2\colon \tilde{M}\to M'$.
Since the space of holomorphic two-forms on a hyperk\"ahler manifold is one-dimensional,
the forms $\pi_1^*\sigma$ and $\pi_2^*\sigma'$ are proportional, we may assume that they are equal,
and denote this two-form by $\tilde{\sigma}$.

Let $\tilde{Z} = \pi_1^{-1}(Z)$, then $Z' = \pi_2(\tilde{Z})$. Every irreducible component
of $Z$ is dominated by an irreducible component of $\tilde{Z}$, so by part (4) of \ref{prop_srank},
we have $\rk_\sigma(Z) \le \rk_{\tilde{\sigma}}(\tilde{Z})$. The inverse inequality follows from
parts (3) and (4) of \ref{prop_srank}, because every irreducible component of $\tilde{Z}$ is mapped
into some irreducible component of $Z$. By an analogous argument applied to $Z'$ and $\tilde{Z}$,
we have $\rk_{\sigma'}(Z') = \rk_{\tilde{\sigma}}(\tilde{Z})$ and this completes the proof.
\end{proof}

\subsection{The symplectic rank of MBM loci and flops}

We use the notion of symplectic rank to study compositions of wall-crossing flops.
This is the key technical step in the proof of the Hu--Yau conjecture.
We start from the following observation.

\hfill

\proposition\label{prop_srank_MBM}
Let $\eta$ be an extremal MBM class on a hyperk\"ahler manifold $M$ and $E_\eta$ its full
MBM locus. Then $\rk_\sigma(E_\eta) = \dim(M) - 2\codim(E_\eta)$.

\hfill

\begin{proof}
By \cite[Proposition 2.2]{_AV:MBM_}, for every irreducible component $E'$ of $E_\eta$
the codimension of $E'$ equals the dimension of the kernel of $\sigma|_{E'}$ at a general
point of $E'$. It follows that $\rk_\sigma(E') = \dim(E') - \codim(E') = \dim(M) - 2\codim(E')$,
which implies the claim, since the symplectic rank of $E_\eta$ is defined by taking the
maximum over all irreducible components, and the codimension is defined by taking the minimum
over all irreducible components.
\end{proof}

\hfill

The following proposition is a strengthening of \ref{prop_WW}.

\hfill

\proposition\label{prop_flop_srank}
Let $\phi\colon M_1 \dra M_2$ be a wall-crossing flop between hyperk\"ahler manifolds of dimension $2n$.
If $M_i$ are non-projective, assume additionally that $b_2(M_i) \ge 6$.
Let $W_1\subset M_1$ and $W_2\subset M_2$ be two closed subvarieties with $\rk_{\sigma_i}(W_i) \le 2n - 6$.
Then there exist two closed subvarieties $Z_1\subset M_1$ and $Z_2\subset M_2$ of codimension
strictly greater than two with $W_i\subset Z_i$, $\rk_{\sigma_i}(Z_i) \le 2n - 6$ and such that $\phi\colon M_1\setminus Z_1\dra M_2\setminus Z_2$ is
Mukai's elementary transformation in codimension two. In particular,
this applies when $W_i$ are MBM loci of codimension strictly greater than two,
or the full transforms of such loci under some bimeromorphic maps.

\hfill

\begin{proof}
First note that by part (2) of \ref{prop_srank} we have $2\codim(W_i) \ge 2n - \rk_{\sigma_i}(W_i)$,
and by our assumptions this implies $\codim(W_i) \ge 3$.

Repeating the proof of \ref{prop_WW} we decompose $\phi$ into the composition
of a contraction $\pi_1\colon M_1\to N$ and the inverse of another contraction $\pi_2\colon M_2\to N$.
We let $W = \pi_1(W_1)\cup \pi_2(W_2)$ and $Z\subset N$ be the union of $W$ and all strata of the Kaledin
stratification of codimension six and higher. Define $Z_i = \pi_i^{-1}(Z)$.
Note that since the symplectic rank of $W$
is at most $2n - 6$, the variety $W$ can not contain any stratum of codimension four
in the Kaledin stratification. Therefore the varieties $Z_1$ and $Z_2$ do not dominate
any stratum of codimension four, so we have $\rk_\sigma(Z_i) \le 2n - 6$.
The latter implies that $\codim(Z_i) \ge 3$, as was explained above, and we conclude the proof as in \ref{prop_WW}.
\end{proof}

\hfill

We are now ready to prove that under certain conditions the compositions of wall-crossing flops
are Hu--Yau transformations. This will give us a basic step in the proof of the Hu--Yau conjecture.

\hfill

\proposition\label{prop_step}
Consider a sequence of wall-crossing flops between $2n$-dimensional hyperk\"ahler manifolds
$$
M_1 \stackrel{\phi_1}{\dra} M_2 \stackrel{\phi_2}{\dra} \ldots \stackrel{\phi_k}{\dra} M_{k+1},
$$
such that exactly one of the MBM loci defining the flops has codimension two. Then the composition
$\phi = \phi_k\circ\ldots\circ\phi_1$ is a Hu--Yau transformation.

\hfill

\begin{proof}
Assume that $\phi_{i_0}$ is a flop in codimension two. Let $W_{i_0 + 1} \subset M_{i_0 + 1}$ be the union of the full transforms
of $E_j$ for $j > i_0$ and let $W_{i_0} \subset M_{i_0}$ be the union of the full transforms
of $E_j$ for $j < i_0$. Then by \ref{prop_srank_inv} the symplectic rank of $W_{i_0}$ and $W_{i_0+1}$ is a most $2n - 6$.
Applying \ref{prop_flop_srank} to the flop $\phi_{i_0}\colon M_{i_0} \dra M_{i_0+1}$ and the two subvarieties
$W_{i_0}$ and $W_{i_0+1}$, we get two subvarieties $Z_{i_0}$ and $Z_{i_0+1}$ of symplectic rank at most $2n-6$
such that $\phi_{i_0}\colon M_{i_0}\setminus Z_{i_0} \dra M_{i_0+1}\setminus Z_{i_0+1}$ is Mukai's elementary transformation
in codimension two, and we have the inclusions $W_{i_0}\subset Z_{i_0}$ and $W_{i_0+1}\subset Z_{i_0+1}$.
We let $Z_1\subset M_1$ be the full transform of $Z_{i_0}$ and $Z_{k+1}\subset M_{k+1}$ be the full transform of $Z_{i_0 + 1}$.
By \ref{prop_srank_inv} the symplectic ranks of $Z_1$ and $Z_{k+1}$ are still a most $2n - 6$, hence they are
of codimension at least three.
The composition $\phi_{i_0-1}\circ\ldots\circ\phi_1\colon M_1\setminus Z_1 \to M_{i_0}\setminus Z_{i_0}$ is a biholomorphic
morphism, because by construction $Z_1$ contains the indeterminacy locus of that map. Analogously, the composition
$\phi_k\circ\ldots\circ\phi_{i_0+1}\colon M_{i_0 + 1}\setminus Z_{i_0 + 1} \to M_{k + 1}\setminus Z_{k + 1}$ is a biholomorphic
morphism. We therefore conclude that $\phi\colon M_1\setminus Z_1 \dra M_{k+1}\setminus Z_{k+1}$ is a composition
of two biholomorphic morphisms and Mukai's elementary transformation in codimension two. Therefore $\phi$ is a Hu--Yau transformation. 
\end{proof}

\section{Proof of the main theorem}\label{sec_main_proof}

We use the results of the previous section to prove \ref{thm_main}.
Given a bimeromorphic map $\varphi\colon M \dra M'$ between compact hyperk\"ahler
manifolds, we apply \ref{prop_decomp} and decompose $\phi$
into a finite sequence of wall-crossing flops $\varphi_i\colon M_i\dra M_{i+1}$,
$i=1,\ldots,k$, such that $M_1 = M$, $M_{k+1} = M'$, so that $\varphi = \varphi_k \circ \ldots \circ \varphi_1$.
We split the sequence of flops $\phi_1, \phi_2, \ldots, \phi_k$ into groups of consecutive flops, each
group containing exactly one flop in codimension two (all the flops in the sequence can not be in codimension greater than two
by the assumptions of the theorem). Applying \ref{prop_step} to
each of the groups of the flops, we see that the composition of the flops forming the group
is a Hu--Yau transformation. Hence the claim of the theorem.

\section{Compositions of flops in codimension two}\label{sec_compositions}

The key step of the proof of \ref{thm_main} is \ref{prop_flop_srank} that allows us
to show that a wall-crossing flop in codimension two is a Mukai elementary transformation
on an open subset with complement of codimension three or higher, even if we require that
the complement should contain a given subvariety (under some suitable conditions on the latter).
The purpose of this section is to show that it is in general not possible to achieve this
for a composition of several flops, i.e. there need not exist an open subset
with complement of codimension higher than two where such map is a composition of successive Mukai's elementary transformations
in codimension two. We will show this on an example.

To construct the example, we start from the K3 surface $S$ that is a quartic hypersurface in $\bbP^3$
containing a line $C$, and generic with this condition. The pencil of planes passing through $C$
defines an elliptic fibration $f\colon S\to \bbP^1$ whose general fibre will be denoted by $F$. We take $M = S^{[3]}$.

Let us recall the description of the MBM classes on $M$ following \cite{_AV:Loci_}, and the
explicit description of the MBM loci that can be found in \cite{HT1}, \cite{HT3} and \cite{_AV:Loci_}.

Let us denote by $\Lambda$ the K3 lattice and let $\Pic(S)\subset \Lambda$ be the Picard lattice of $S$.
The latter is spanned by the classes of the curves $C$ and $F$ and has the intersection matrix
$$
\Pic(S) = \begin{pmatrix}
-2 & 3\\
3 & 0
\end{pmatrix}
$$
Let $\Lambda'$ be the lattice $H^2(M, \bbZ)$ with the BBF form. Then $\Lambda'\simeq \Lambda\oplus\bbZ\delta$,
where the direct sum is orthogonal and $\delta = [\Delta]/2$ with $\Delta$ the exceptional divisor of the contraction $S^{[3]}\to S^{(3)}$.
We have $\delta^2 = -4$. Consequently, we also have $\Pic(M) \simeq \Pic(S)\oplus \bbZ\delta$.
The lattice $\Lambda$ is unimodular and the discriminant group of the lattice $\Lambda'$ is clearly isomorphic to $\bbZ/4\bbZ$.

According to \cite{_AV:Loci_} the MBM classes on $M$ form five orbits under the action of the diffeomorphism
group of $M$: three orbits consisting of divisorial classes, one orbit of codimension two classes and one orbit of
codimension three classes. The orbits are distinguished by the square of the class $x\in \Lambda'$ and by the
image of $x/d(x)$ in the discriminant group, where $d(x)$ is the divisibility of $x$ as an element of the dual
lattice $(\Lambda')^\vee$. More precisely, $x_1$ and $x_2$ are in the same orbit if and only if $x_1^2=x_2^2$
and $x_1/d(x_1) \equiv \pm x_2/d(x_2)\, (\mathrm{mod}\,\, \Lambda')$.

The representatives of the orbits in our situation are the following:
\begin{enumerate}
\item The divisorial class $\delta$; we have $\delta^2 = -4$ and $\delta/d(\delta) = \delta/4$ is the class of the
ruling on the divisor $\Delta = E_\delta$;
\item The divisorial class $\alpha = [C]$; we have $\alpha^2 = -2$ and $\delta(\alpha) = 1$. The divisor $E_\alpha$
consists of subschemes that intersect $C$;
\item The divisorial class $\epsilon = 2[F] - \delta$; we have $\epsilon/d(\epsilon) = [F] - \delta /2$ is the class
of the ruling on the divisor $E_\epsilon$. The latter divisor is the closure of the image
of the natural map $\mathrm{Hilb}^2(S/\bbP^1)\times S\dra S^{[3]}$, where $\mathrm{Hilb}^2(S/\bbP^1)$
is the relative Hilbert scheme of the elliptic fibration $f\colon S\to \bbP^1$;
\item Let us describe three classes of codimension two.
\begin{enumerate}
\item The class $\beta = 4[C] - \delta$ with $\beta^2 = -36$ and $d(\beta) = 4$. The MBM locus $E_\beta$
is the closure of the image of the natural map $C^{(2)}\times S\dra S^{[3]}$;
\item The class $\eta = 4[F] - 3\delta$ with $\eta^2 = -36$ and $d(\eta) = 4$. The MBM locus $E_\eta$
is the image of the natural embedding $\mathrm{Hilb}^3(S/\bbP^1)\hrarr S^{[3]}$;
\item The class $\zeta = 4[C] + 4[F] - 5\delta$ with $\zeta^2 = -36$ and $d(\zeta) = 4$. The MBM locus
$E_\zeta$ consists of subschemes of $S$ that lie on a line in $\bbP^3$; see Remark \ref{rem_MBM} below;
\end{enumerate}
\item The codimension three class $\gamma = 2[C] - \delta$; we have $\gamma^2 = -12$ and $\delta(\gamma) = 2$.
The MBM locus $E_\gamma$ is the set of subschemes supported on $C$, isomorphic to $C^{(3)}\simeq \bbP^3$.
\end{enumerate}

\begin{rem}\label{rem_MBM}
	Descriptions of all MBM classes and the corresponding rational curves from the above list have already
	appeared in the literature, see \cite[section 4]{_AV:Loci_} and references therein. One exception is,
	possibly, the class $\zeta$. Let us indicate how one can find the expression of this class given above.
	A rational curve in the MBM locus $E_\zeta$ can be described as follows. Let us fix a
	line $L_0\subset \bbP^3$ such that $L_0 \cap S = \{p, q, r, s\}$, where $p \in C$ and $\{q, r, s\}$ are
	three other distinct point on $S$. Let $P_0\subset \bbP^3$ be a plane passing through the line $L_0$
	and not containing $C$.
	Consider all lines $L\subset P_0$ passing through the point $q$. For a generic such line the intersection
	$S\cap L$ consists of the point $q$ and three other points on $S$, the latter defining
	a point of $S^{[3]}$. This gives a rational curve $C_\zeta$ in $S^{[3]}$. To find its cohomology class,
	we can compute its intersection indices with the divisor classes $[C]$, $[F]$ and $\epsilon$.
	We have chosen the line $L_0$ in a way to make sure that the points of intersection are smooth points of the
	divisors and the intersections are transversal.
	The divisor with the class $[C]$ is $E_\alpha$, it consists of all subschemes of $S$ intersecting the
	curve $C$. By construction, there is only one such subscheme $\{p, r, s\}$ on the curve $C_\zeta$, hence
	$C_\zeta \cdot [C] = 1$. Analogously, to compute $C_\zeta\cdot[F]$, consider a generic plane $P_1\subset \bbP^3$
	containing $C$. Let $F_1$ be the irreducible component of $P_1\cap S$ of degree three. The class
	$[F]$ on $S^{[3]}$ is represented by the divisor consisting of subschemes intersecting $F_1$.
	Note that $F_1\cap P_0$ are three points on a line in $P_0$. For a generic choice of $P_1$ the latter
	line is different from $L_0$, and therefore $C_\zeta$ has three distinct points of intersection with
	the above divisor, therefore $C_\zeta\cdot[F] = 3$. To compute $C_\zeta\cdot \epsilon$, note that the
	unique point of intersection of $C_\zeta$ and $E_\epsilon$ is the subscheme $\{p, r, s\}$ corresponding
	to the line $L_0$, hence $C_\zeta\cdot \epsilon = 1$. Using these intersection indices, we find that
	$[C_\zeta] = [C] + [F] - \frac{5}{4}\delta$, hence the expression for the MBM class $\zeta$ given
	in the list above.
\end{rem}

\begin{figure}[h]
\centering
\includegraphics[width=10cm]{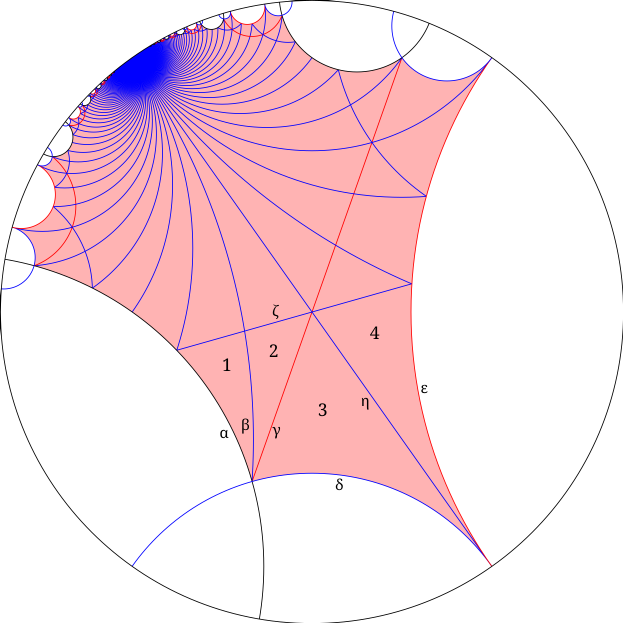}
\caption{The birational ample cone of $S^{[3]}$}\label{fig_1}
\end{figure}

In the above listing we have slightly abused the notation for the MBM loci: not all of the MBM classes
on the list are extremal, so not all of the loci described above are contractible on $S^{[3]}$. In order
to make them contractible one may need to pass to a different birational model. The birational
ample cone of $S^{[3]}$ (i.e. the intersection of the birational K\"ahler cone with the subspace
spanned by the Picard lattice) and its wall-and-chamber structure is depicted in Figure \ref{fig_1}.

In Figure \ref{fig_1} we have the projectivisation of the positive cone, which is the hyperbolic plane in our case,
and the walls described above. We will consider the K\"ahler chambers numbered 1, 2, 3 and 4, and the corresponding
birational models $M_i$ for $i=1,2,3,4$. The model $M_3$ is the Hilbert scheme $S^{[3]}$.
The walls are colored depending on the image of $x/d(x)$ in the discriminant group for an MBM class $x$:
the black walls correspond to $0 \, (\mathrm{mod}\,\, 4)$, the blue walls to $\pm 1 \, (\mathrm{mod}\,\, 4)$
and the red walls to $2 \, (\mathrm{mod}\,\, 4)$.

The manifold $S^{[3]}$ admits a Lagrangian fibration induced by the elliptic fibration $f\colon S\to \bbP^1$,
the nef isotropic class defining the fibration being equal to $[F]$ (the lower right cusp in Figure \ref{fig_1}).
It also admits a rational Lagrangian fibration defined by the isotropic class $[C] + [F] - \delta$ (the point
of concentration of the walls in the upper left part of Figure \ref{fig_1}). This fibration is due to the fact
that $S^{[3]}$ is birational to the Beauville--Mukai system with the Mukai vector $(0, [C]+[F], 1)$; it can be described
explicitly as the map $S^{[3]}\dra \bbP^3$ sending a triple of points on $S$ to the plane in $\bbP^3$
spanned by these points. One can see that the indeterminacy locus of this map is the MBM locus $E_\zeta$,
and one passes to a model where the fibration becomes regular by making two flops: in $E_\gamma \simeq \bbP^3$
(which is the singular locus of $E_\zeta$) and then in $E_\zeta$. The Beauville--Mukai system in our example admits an infinite
group of birational automorphisms, hence the infinite family of walls passing
through the isotropic class $[C] + [F] - \delta$.

We will consider the following sequence of flops:
$$
M_1\stackrel{\phi_\beta}{\dra} M_2 \stackrel{\phi_\gamma}{\dra} M_3 \stackrel{\phi_\eta}{\dra} M_4,
$$
denoting by $\phi$ their composition. Both $\phi_\beta$ and $\phi_\eta$ are wall-crossing flops in codimension two,
while $\phi_\gamma$ is a wall-crossing flop in codimension three.

To understand the geometry of the map $\phi$, we need to study the intersections of the MBM loci
that define the flops. First, consider the subvariety $E_\beta$. It is non-normal and
birational to $\bbP^2 \times S \simeq C^{(2)}\times S$. The singular locus of $E_\beta$ is $E_\gamma \simeq C^{(3)} \simeq \bbP^3$. Denote by $\tilde{E}_\beta$
the normalization of $E_\beta$, then $\tilde{E}_\beta$ is isomorphic to the blow-up of
$C^{(2)}\times S$ in the incidence subvariety $Z\subset C^{(2)} \times C\subset C^{(2)}\times S$,
where $((x,y), z)\in Z$ iff $x=z$ or $y=z$. The proper transform of $C^{(2)}\times C$ under the
blow-up is a divisor $D\subset \tilde{E}_\beta$. The proper transform of $E_\beta$ on the model $M_2$,
that we will denote by $E'_\beta$, is obtained by blowing-down the divisor $D$. Then variety $E'_\beta$
is a $\bbP^2$-bundle over $S$ and it is the center of the flop $\phi_\beta^{-1}$.

The center of the flop $\phi_\eta$ is $E_\eta\simeq \mathrm{Hilb}^3(S/\bbP^1)$,
which is a $\bbP^2$-bundle over another K3 surface $S' = \Pic^3(S/\bbP^1)$, the latter being the
relative Picard group of the elliptic fibration $f\colon S\to \bbP^1$. For $t\in \bbP^1$ we will
denote by $F_t$ the fiber of $f$ over $t$. We will implicitly fix an origin on $F_t$ making
it an elliptic curve and also identifying the fiber of $S'$ over $t$ with $F_t$.

Let $E_{\beta\eta}=E_\beta\cap E_\eta\subset M_3$. Then we have a morphism $E_{\beta\eta}\to S'$
that is the composition of the embedding into $E_\eta$ and the projection onto $S'$. Analogously,
we get a rational map $E_{\beta\eta} \dra S$ given by the embedding into $E_\beta$ and
the rational projection onto $S$. We will denote by $\Sigma\subset S\times S'$ the closure of the image
of the product of the two maps above. Note that $\Sigma$ is fibered over $\bbP^1$ and
denote by $\Sigma_t = \Sigma\cap (F_t\times F_t)$ a fiber.

\hfill

\lemma\label{lem_cor}
The subvariety $\Sigma\subset S\times S'$ is dominant and finite of degree three both over $S$ and $S'$.
If the fiber $F_t\subset S$ intersects the curve $C$ in three distinct points $e_0$, $e_1$ and $e_2$,
then $\Sigma_t\subset F_t\times F_t$ consists of three connected components that are the translations
of the diagonal by $(0, e_0 + e_1)$, $(0, e_1 + e_2)$ and $(0, e_2 + e_0)$.

\hfill

\begin{proof}
Let us fix a point $t\in \bbP^1$ so that $F_t$ satisfies the conditions of the lemma. Let $x\in F_t$
be a point different from $e_0$, $e_1$, $e_2$. Then the fiber of $E_\beta$ over $x$ is $\bbP^2\simeq C^{(2)}\times \{x\}$.
It intersects $E_\eta$ in three points: $((e_0, e_1), x)$, $((e_1, e_2), x)$ and $((e_2, e_0), x)$
that are mapped by the projection $E_\eta\to S'$ to the points $e_0 + e_1 + x$, $e_1 + e_2 + x$
and $e_2 + e_0 + x$ respectively. This proves the claim about $\Sigma_t$, because the latter is closed in $F_t\times F_t$.
Hence $\Sigma_t$ is dominant and finite of degree three over both factors, and since this holds for a generic $t\in \bbP^1$,
we also get the first claim of the lemma.
\end{proof}

\hfill

Denote by $p$ and $q$ the projections from $\Sigma$ to the two factors of $S\times S'$.
Given $x, y\in S$ we will write $x\sim y$ if $q(p^{-1}(x))\cap q(p^{-1}(y)) \neq \emptyset$.
Generating an equivalence relation on the points of $S$ by $\sim$, we will call the corresponding
equivalence classes $\Sigma$-orbits.

\hfill

\lemma\label{lem_equiv}
Let $x\in F_t$. Then $x\sim x + e_1 - e_0$, $x\sim x + e_2 - e_1$ and $x \sim x + e_0 - e_2$.
Whenever at least one of $e_1 - e_0$, $e_2 - e_1$ or $e_0 - e_2$ is non-torsion,
the orbit of $x$ is infinite, hence Zariski dense in $F_t$.

\hfill

\begin{proof}
From \ref{lem_cor} we see that $$q(p^{-1}(x)) = \{x + e_0 + e_1, x + e_1 + e_2, x + e_2 + e_0\},$$
and $$q(p^{-1}(x + e_1 - e_0)) = \{x + 2e_1, x + 2e_1 - e_0 + e_2, x + e_1 + e_2\},$$ which
shows that $x \sim x + e_1 - e_0$. Analogously for the other two points. The last claim
of the lemma is obvious.
\end{proof}

\hfill

\lemma\label{lem_nontors}
For a very general $t\in \bbP^1$ the points $e_1 - e_0$, $e_2 - e_1$ and $e_0 - e_2$
are non-torsion in $F_t$.

\hfill

\begin{proof}
Let $U\subset \bbP^1$ be the complement of the discriminant of the elliptic fibration $f$
and let $g\colon J\to U$ be the Jacobian fibration corresponding to $f|_U$. Let $J^t\subset J$
be the subset of points that are torsion in their respective fibers. Then $J^t$ is a countable
union of closed subvarieties of $J$, each of them being an unramified finite cover of $U$.

For a generic choice of the surface $S$ the curve $C$ is ramified over $\bbP^1$ so that at least one of the
ramification points is contained in a smooth fiber of the elliptic fibration $f$.
The pairwise differences of the points $e_0$, $e_1$ and $e_2$ (that are the intersection points
of $F_t$ with $C$) define a curve $C'\subset J$ that is a ramified covering of $U$.
Hence $C'$ can not be contained in $J^t$, therefore $J^t\cap C'$ and its image in $U$ are
countable sets. This proves the claim.
\end{proof}

\hfill

As the final step of our discussion let us show that the birational map $\phi$ can not be represented
as a composition of Mukai flops in codimension two, even after removing some subvarieties of higher
codimension from $M_1$ and $M_4$, see \ref{prop_no_open} below. Before proving this claim, let us
introduce some notation. Let $Z_\beta$ be the center of the flop $\phi_\beta$,
$Z_\eta$ be the proper transforms of $E_\eta$ on the model $M_1$,
and $Z_\gamma$ be the proper transform of the center of the flop $\phi_\gamma$ on the model $M_1$. By construction,
$Z_\beta$, $Z_\gamma$ and $Z_\eta$ are the irreducible components of the indeterminacy locus
of the map $\phi$. Note that $Z_\gamma$ is of codimension three, being birational to $\bbP^3$,
and the other two components are of codimension two. The subvarieties $Z_\beta$, $Z_\gamma$ and $Z_\eta$
are swept out by rational curves representing the corresponding MBM classes.

Assume that there exist two open subsets $M_1^\circ\subset M_1$ and $M_4^\circ\subset M_4$
such that $\phi|_{M_1^\circ} = \phi_m\circ\ldots \circ\phi_0$,
where $\phi_i\colon U_i\dra U_{i+1}$ are Mukai's elementary transformations in codimension two,
$U_0 = M_1^\circ$ and $U_{m+1} = M_4^\circ$. Every $\phi_i$ is the composition
of a proper contraction $U_i\to V_i$ with the exceptional set $E_i$ and the inverse
of a proper contraction $U_{i+1}\to V_i$ with the exceptional set $E'_{i+1}$.
Both $E_i$ and $E'_{i+1}$ are smooth $\bbP^2$-bundles over the singular locus of $V_i$.
We may also assume that $E_i$ and $E'_{i+1}$ are irreducible.

Let us denote by $Z_i$ the proper transforms of the varieties $E_i$ on the model $M_1^\circ$.
By this we mean the following. Given an irreducible subvariety of codimension two $Z\subset U_{i+1}$,
we have two possibilities. Either $Z$ is not contained in $E'_{i+1}$, and then we take its
proper preimage under $\phi_i$ in the usual sense: the closure in $U_i$ of the preimage
of $Z\setminus E_{i+1}'$ under $\phi_i$. Or $Z$ coincides with $E'_{i+1}$, and
in this case by the proper preimage of $Z$ under $\phi_i$ we will mean the subvariety $E_i$.
Iterating this procedure over the maps $\phi_i$, we get a proper preimage of $Z$ in $M_1^\circ$.

The proper preimages $Z_i$ are codimension two irreducible subvarieties of $M_1^\circ$
that are swept out by rational curves $C_i$, the latter being proper preimages of the
lines in the fibers of the $\bbP^2$-bundles $E_i$. Let us denote by $\overline{Z}_i$ the
closure of $Z_i$ in $M_1$. Now we can prove the following statement.

\hfill
 
\proposition\label{prop_no_open}
There do not exist open subsets $M_1^\circ\subset M_1$ and $M_4^\circ\subset M_4$ with complements
of codimension strictly greater than two and such that $\phi\colon M_1^\circ \dra M_4^\circ$
is a composition of Mukai's elementary transformations in codimension two.

\hfill

\begin{proof} Assume the contrary and use the notation introduced above. To arrive at a contradiction,
we will make two steps. First, we will show that the subvariety $Z_\gamma$ is contained in the complement
of $M_1^\circ$. Second, we will use \ref{lem_cor}, \ref{lem_equiv} and \ref{lem_nontors} to show that
$Z_\beta$ and $Z_\eta$ also have to be in the complement of $M_1^\circ$ contradicting the assumption on
the codimension.

{\it Step 1.} Note that the indeterminacy locus of the map $\phi|_{M_1^\circ}$ has to be contained in the union
of the subvarieties $Z_i$ constructed before the statement of the proposition.
Assume that $Z^\circ_\gamma = Z_\gamma\cap M_1^\circ\neq \emptyset$. Then $Z^\circ_\gamma$ has
to be contained in one of $Z_i$. Let $i_0$ be the smallest index such that $Z^\circ_\gamma\subset Z_{i_0}$,
so that $Z_\gamma\subset\overline{Z}_{i_0}$. Let $Z_\gamma'\subset U_{i_0}$ be the proper transform
of $Z_\gamma^\circ$ under the sequence of flops $\phi_0,\ldots,\phi_{i_0-1}$. By construction, $Z_\gamma'\subset E_{i_0}$.

We claim that $Z_\gamma'$ does not dominate the base of the $\bbP^2$-bundle $E_{i_0}$. To prove the claim, let
$\widetilde{Z}_{i_0}$ be the resolution of singularities of $\overline{Z}_{i_0}$. Observe that the singular
locus of $Z_{i_0}$ is contained in the union of $Z_i$ with $i < i_0$. Therefore, our choice of $i_0$
guarantees that $Z_\gamma^\circ$ is not contained in the singular locus of $Z_{i_0}$. We can then
consider the proper preimage $\widetilde{Z}_\gamma$ of $Z_\gamma$ in $\widetilde{Z}_{i_0}$,
and the variety $\widetilde{Z}_\gamma$ is then birational to $Z_\gamma$. Since $Z_\gamma$ is rational,
the variety $\widetilde{Z}_\gamma$ can not dominate the base of the MRC-fibration of $\widetilde{Z}_{i_0}$.
The latter is birational to the $\bbP^2$-bundle $E_{i_0}$, hence the claim that $Z_\gamma'$ does not dominate
the base of that $\bbP^2$-bundle.

The previous paragraph implies that the variety $Z_\gamma'$ contains one of the fibers of the $\bbP^2$-bundle $E_{i_0}$.
It follows that one of the rational curves $C_{i_0}$ defined above
may be chosen to lie inside the variety $Z_\gamma$. Then the image of $C_{i_0}$ on the model $M_3$ is
contained in the MBM locus $E_\gamma\simeq \bbP^3$, therefore all deformations of $C_{i_0}$ stay inside $E_\gamma$.
On the other hand, by construction, the deformations of $C_{i_0}$ sweep out a subvariety of codimension two.
This contradiction completes the first step of the proof.

{\it Step 2.} Denote by $\pi_1\colon M_1\to N_1$ and $\pi_4\colon M_4\to N_4$ the extremal contractions
of the MBM loci that define the wall-crossing flops $\phi_\beta$ and $\phi_\eta$.
Note that the exceptional set of $\pi_1$ is the variety $Z_\beta$, and since it has codimension two,
it has to coincide with $\overline{Z}_i$ for some $i$. The variety $Z_i = \overline{Z}_i\cap M_1^\circ = Z_\beta\cap M_1^\circ$ is a proper
$\bbP^2$-bundle, hence $M_1^\circ$ has to be proper over its image in $N_1$.
Analogously, $M_4^\circ$ has to be proper over its image in $N_4$.

The singular locus of $N_1$ is the K3 surface $S$. Assume that
$x\in S$ is a point in the complement of $\pi_1(M_1^\circ)$. Our definition of
the correspondence $\Sigma\subset S\times S'$ from \ref{lem_cor} implies that the
entire $\Sigma$-orbit of $x$ has to lie in the complement of $\pi_1(M_1^\circ)$.
By \ref{lem_equiv} and \ref{lem_nontors} such orbit will be dense in the fiber $F_t$ of the elliptic
fibration $f$ if $t\in \bbP^1$ is very general.

Now recall that by the first step of the proof $Z_\gamma$ has to be contained in the complement of $M_1^\circ$.
Therefore the curve $C\subset S \subset N_1$ is in the complement of $\pi_1(M_1^\circ)$.
The curve $C$ intersects all fibers of the elliptic fibration $f$, so by the previous
paragraph the complement of $\pi_1(M_1^\circ)$ contains a very general fiber $F_t$.
Therefore it has to contain the whose surface $S$, because very general fibers form a dense subset of $S$.
We conclude that the entire MBM locus $Z_\beta$ has to be in the complement of $M_1^\circ$.
But $Z_\beta$ has codimension two in $M_1$, hence the claim of the proposition.
\end{proof}

\hfill

{\bf Acknowledgements:} 
We are grateful to Yulia Gorginyan and Dmitry Kaledin for interesting discussions of the subject.

\hfill

{

}

{\small
\noindent {\sc Ekaterina Amerik\\
{\sc Laboratory of Algebraic Geometry,\\
National Research University HSE,\\
Department of Mathematics, 6 Usacheva Str. Moscow, Russia,}\\
\tt  Ekaterina.Amerik@gmail.com}, also: \\
{\sc Universit\'e Paris-11,\\
Laboratoire de Math\'ematiques,\\
Campus d'Orsay, B\^atiment 425, 91405 Orsay, France}

\hfill

\noindent {\sc Andrey Soldatenkov\\
{\tt aosoldatenkov@gmail.com}\\
{\sc Universidade Estadual de Campinas\\
Departamento de Matem\'atica - IMECC\\
Rua S\'ergio Buarque de Holanda, 651\\
13083-859, Campinas - SP, Brasil}}

\hfill

\noindent {\sc Misha Verbitsky\\
{\sc Instituto Nacional de Matem\'atica Pura e
              Aplicada (IMPA) \\ Estrada Dona Castorina, 110\\
Jardim Bot\^anico, CEP 22460-320\\
Rio de Janeiro, RJ - Brasil\\
\tt  verbit@impa.br }
}}

\end{document}